\documentclass[a4paper,10pt]{amsart}

\usepackage[latin1]{inputenc}
\usepackage{amsmath}
\usepackage{amsfonts}
\usepackage{amssymb}
\usepackage[francais,english]{babel}
\usepackage{dsfont}
\usepackage[all]{xy}
\usepackage{stmaryrd}

\newcommand{\al}{\alpha}
\newcommand{\eps}{\varepsilon}

\newcommand{\fl}{\rightarrow}
\newcommand{\gfl}{\longrightarrow}

\newcommand{\dr}{\ar@{->}[r]}
\newcommand{\drp}{\ar@{-->}[r]}
\newcommand{\ddr}{\ar@{->}[rr]}
\newcommand{\ha}{\ar@{->}[u]}
\newcommand{\hae}{\ar@{->>}[u]}
\newcommand{\hap}{\ar@{-->}[u]}
\newcommand{\ham}{\ar@{^{(}->}[u]}
\newcommand{\hham}{\ar@{^{(}->}[uu]}
\newcommand{\bdr}{\ar@{->}[dr]}
\newcommand{\bbdr}{\ar@{->}[ddr]}
\newcommand{\bddr}{\ar@{->}[drr]}
\newcommand{\hdr}{\ar@{->}[ur]}
\newcommand{\hdrp}{\ar@{-->}[ur]}
\newcommand{\bas}{\ar@{->}[d]}
\newcommand{\basp}{\ar@{-->}[d]}
\newcommand{\bg}{\ar@{->}[dl]}
\newcommand{\bgp}{\ar@{-->}[dl]}

\newcommand{\cat}{\mathcal{C}}
\newcommand{\ec}{\mathcal{E}}
\newcommand{\mc}{\mathcal{M}}
\newcommand{\pc}{\mathcal{P}}

\newcommand{\tc}{\mathcal{T}}

\newcommand{\ms}{\underline{\mathcal{M}}}
\newcommand{\aq}{\mathcal{A}_Q}
\newcommand{\aqt}{\mathcal{A}_{Q_T}}
\newcommand{\aqtp}{\mathcal{A}_{Q_{T'}}}

\newcommand{\modm}{\operatorname{Mod}\mathcal{M}}
\newcommand{\modms}{\operatorname{mod}\underline{\mathcal{M}}}
\newcommand{\perm}{\operatorname{per}\mathcal{M}}
\newcommand{\perms}{\operatorname{per}_{\underline{\mathcal{M}}}\mathcal{M}}

\newcommand{\z}{\mathbb{Z}}

\newcommand{\susp}{\Sigma}
\newcommand{\susm}{\Sigma^{-1}}
\newcommand{\ke}{\operatorname{Ker}}
\newcommand{\im}{\operatorname{Im}}
\newcommand{\coke}{\operatorname{Coker}}

\newcommand{\ind}{\operatorname{ind}}
\newcommand{\coi}{\operatorname{coind}}
\newcommand{\gre}{\operatorname{Gr}_e}
\newcommand{\grf}{\operatorname{Gr}_f}
\newcommand{\grg}{\operatorname{Gr}_g}
\newcommand{\ko}{\text{\rm K}_0(\operatorname{mod}B)}
\newcommand{\kosp}{\text{\rm K}_0^{\text{sp}}(\operatorname{mod}B)}
\newcommand{\kop}{\text{\rm K}_0(\operatorname{proj}B)}
\newcommand{\modb}{\operatorname{mod}B}

\newcommand{\ext}{\operatorname{Ext}}
\newcommand{\homph}{\operatorname{Hom}}
\newcommand{\homb}{\operatorname{Hom}_B}

\newtheorem{theo}{Theorem}[section]
\newtheorem{prop}[theo]{Proposition}
\newtheorem{lem}[theo]{Lemma}
\newtheorem{cor}[theo]{Corollary}
\newtheorem{defi}[theo]{Definition}

\begin{document}

\title{Cluster characters for 2-Calabi--Yau triangulated categories}
\author{Yann Palu}
\address{Universit{\'e} Paris 7 - Denis Diderot, 
UMR 7586 du CNRS, case 7012, 2 place Jussieu, 
75251 Paris Cedex 05, France.}

\email{
\begin{minipage}[t]{5cm}
palu@math.jussieu.fr
\end{minipage}
}

\begin{abstract}
Starting from an arbitrary
cluster-tilting object $T$ in a
2-Calabi--Yau category over an algebraically
closed field, as in the setting of Keller and Reiten, 
we define, for each object $L$, a fraction
$X(T,L)$ using a formula proposed by Caldero and Keller.
We show that the map taking $L$ to $X(T,L)$ 
is a cluster character, i.e. that it satisfies a 
certain multiplication formula. We deduce that it
induces a bijection, in the finite 
and the acyclic case, between the indecomposable 
rigid objects of the cluster category and
the cluster variables, which confirms a 
conjecture of Caldero and Keller.
\end{abstract}

\maketitle

\section*{Introduction}

Cluster algebras were invented and studied by S. Fomin and A. Zelevinsky in 
\cite{FZ1}, \cite{FZ2}, \cite{FZcoef} 
and in collaboration with A. Berenstein in~\cite{BFZ}.
They are commutative algebras endowed with a distinguished
set of generators called the cluster variables.
These generators are gathered into overlapping sets 
of fixed finite cardinality, called clusters,
which are defined recursively from an initial one via an operation called mutation.
A cluster algebra is said to be of finite type if
it only has a finite number of cluster variables.
The finite type cluster algebras were classified in \cite{FZ2}. 

It was recognized in \cite{MRZ} that the combinatorics of cluster mutation
are closely related to those of tilting theory
in the representation theory of quivers and finite dimensional algebras.
This discovery was the main motivation for the
invention of cluster categories (in \cite{CCS} for the $A_n$-case
and in \cite{BMRRT} for the general case). 
These are certain triangulated categories \cite{Ktri}
which, in many cases, allow one to `categorify' cluster algebras:
In the categorical setting, the cluster-tilting objects play the role 
of the clusters, and their indecomposable direct summands 
the one of the cluster variables. 

In \cite{GLS1}, \cite{GLS2}, \cite{GLS3}, the authors study another setting for
the categorification of cluster algebras: The module categories 
of preprojective algebras of Dynkin type. They succeed in 
categorifying a different class of cluster algebras,
which also contains many cluster algebras of infinite type.

Both cluster categories and module categories of preprojective algebras
of Dynkin type are 2-Calabi--Yau categories in the sense that
we have bifunctorial isomorphisms
$$
\ext^1(X,Y) \simeq D\ext^1(Y,X),
$$
which are highly relevant in establishing the link with
cluster algebras. This motivates the study of 
more general 2-Calabi--Yau categories in
\cite{KR1}, \cite{KR2}, \cite{Tab}, \cite{KZ}, \cite{IR}, \cite{IY}, \cite{BIRS}.
In order to show that a given 2-Calabi--Yau category
`categorifies' a given cluster algebra, a crucial point is 
\begin{enumerate}
 \item[a)] to construct an explicit map from the set of 
	indecomposable factors of cluster-tilting objects 
	to the set of cluster variables, and
 \item[b)] to show that it is bijective.
\end{enumerate}
Such a map was constructed for module categories of 
preprojective algebras of Dynkin type in \cite{GLS1}
using Lusztig's work \cite{L}. 
For cluster categories, it was defined by P. Caldero
and F. Chapoton in \cite{CC}.
More generally, for each object $M$ of 
the cluster category, they defined a fraction $X_M$
in $\mathbb{Q}(x_1,\ldots,x_n)$. The bijectivity
property of the Caldero--Chapoton map was proved
in \cite{CC} for finite type and in \cite{CK2}, 
cf. also~\cite{BCKMRT}, for acyclic type.

A crucial property of the Caldero--Chapoton map is the following.
For any pair of indecomposable objects $L$ and $M$ of $\cat$ whose 
extension space $\cat(L,\susp M)$ is one-dimensional, we have
$$
X_L X_M = X_B + X_{B'} ,
$$
where $\susp$ denotes the suspension in $\cat$ and 
where $B$ and $B'$ are the middle terms of `the' two
non-split triangles with outer terms $L$ and $M$. 
We define, in defintion~\ref{def:clchar} a cluster character 
to be a map satisfying this multiplication formula.

This property has been proved in \cite{CK1} in the finite case,
in \cite{GLS3} for the analogue of the Caldero--Chapoton map 
in the preprojective case, 
and in \cite{CK2} in the acyclic case.

The main result of this article is a generalisation of this
multiplication formula. Starting from an arbitrary
cluster-tilting object $T$ and an arbitrary
2-Calabi--Yau category $\cat$ over an algebraically
closed field (as in the setting of \cite{KR1}),
we define, for each object $L$ of $\cat$, a fraction
$X^T_L$ using a formula proposed in \cite[6.1]{CK1}.
We show that the map $L \mapsto X^T_L$ is a 
cluster character. We deduce that it
has the bijectivity property in the finite 
and the acyclic case, which confirms
conjecture~2 of \cite{CK1}.
Here, it yields a new way of expressing cluster variables
as Laurent polynomials in the variables of a 
fixed cluster.
Our theorem also applies to stable categories of
preprojective algebras of Dynkin type and their 
Calabi--Yau reductions studied in \cite{GLS4} and \cite{BIRS}.

Let $k$ be an algebraically closed field, and let $\cat$ 
be a 2-Calabi--Yau Hom-finite triangulated $k$-category 
with a cluster-tilting object $T$ (see section \ref{s:main}).

The article is organised as follows:
In the first section, the notations are given and the main result is stated.
In the next two sections, we investigate the exponents
appearing in the definition of $X^T_L$. In section \ref{s:ind},
we define the index and the coindex of an object of $\cat$
and show how they are related to the exponents.
Section~\ref{s:abf} is devoted to the study of the antisymmetric
bilinear form
$$
\langle\;,\;\rangle_a
$$
on {\sf mod}$\,$End$_{\cat}T$.
We show that this form descends to the 
Grothendieck group K$_0({\sf mod}\,$End$_{\cat}T)$, 
confirming conjecture 1 of \cite[6.1]{CK1}.
In section \ref{s:dicho}, we prove that the same phenomenon of dichotomy
as in \cite[section 3]{CK2} (see also \cite{GLS3}) still holds in our setting. The results 
of the first sections are used in section \ref{s:mf} 
to prove the multiplication formula. 
We draw some consequences in section \ref{ss:csq}. 
Two examples are given in section~\ref{ss:examples}.

\section*{Acknowledgements}
This article is part of my PhD thesis, under the supervision of Professor
B. Keller. I would like to thank him deeply for introducing me to the subject. 
I would also like to thank P. Caldero for his comments, and J. Schr\"oer for suggesting 
lemma~\ref{Schroer} to me.

\tableofcontents

\section{Main result}\label{s:main}

Let $k$ be an algebraically closed field, and let $\cat$ be a 
$k$-linear triangulated category with split idempotents. 
Denote by $\susp$ its suspension functor. 
Assume moreover that the category $\cat$ 
\begin{itemize}
\item[a)] is Hom-finite: For
any two objects $X$ and $Y$ in $\cat$, the space of morphisms 
$\cat(X,Y)$ is finite-dimensional,
\item[b)] is 2-Calabi--Yau: There exist bifunctorial isomorphisms
$$
\cat(X,\susp Y) \simeq D\cat(Y,\susp X),
$$
where $D$ denotes the duality functor ${\sf Hom}_k(?,k)$, and
\item[c)] admits a cluster-tilting object $T$, which means that
\begin{itemize}
\item[i)] $\cat(T,\susp T) = 0$ and
\item[ii)] for any $X$ in $\cat$, if $\cat(X,\susp T) = 0$, then $X$ belongs
to the full subcategory add$\,T$ formed by the direct summands of sums 
of copies of $T$.
\end{itemize}
\end{itemize}

For two objects $X$ and $Y$ of $\cat$, 
we often write $(X,Y)$ for the space of morphisms $\cat(X,Y)$ 
and we denote its dimension by $[X,Y]$. 
Similarly, we write $^1(X,Y)$ for $\cat(X,\susp Y)$ 
and $^1[X,Y]$ for its dimension. 
Let $B$ be the endomorphism algebra of $T$ in $\cat$, 
and let $\modb$ be the category of finite-dimensional 
right $B$-modules.
As shown in \cite{BMR1}, cf. also \cite{KR1}, the  functor 

$$
F : \cat  \gfl  \modb \;\text{ , } \; X  \longmapsto \cat(T,X),
$$

\noindent induces an equivalence of categories

$$
\cat/(\susp T) \stackrel{\simeq}{\gfl} \modb,
$$

\noindent where $(\susp T)$ denotes the ideal of morphisms of $\cat$ 
which factor through a direct sum of copies of $\susp T$. 

The following useful proposition is proved in \cite{KR1} and \cite{KZ}:

\begin{prop}
Let $ X \stackrel{f}{\fl} Y \stackrel{g}{\fl} Z \fl \susp X$ be a triangle in $\cat$. Then
\begin{itemize}
\item[.] The morphism $g$ induces a monomorphism in $\modb$ if and only if $f~\!\!\in~\!\!(\susp T)$.
\item[.] The morphism $f$ induces an epimorphism in $\modb$ if and only if $g\in(\susp T)$
\end{itemize}
Moreover, if $X$ has no direct summands in add$\,\susp T$, then
$FX$ is projective (resp. injective) if and only if $X$ lies in {\rm add}$\,(T)$ (resp. in 
{\rm add}$\,(\susp^2T)$ ). 
\end{prop}

\begin{defi}\label{def:clchar}
 A cluster character on $\cat$ with values 
in a commutative ring $A$ is a map 
$$
\chi : \text{ obj}(\cat) \gfl A
$$
such that 
\begin{itemize}
 \item[.] for all isomorphic objects $L$ and $M$, we have $\chi(L) = \chi(M)$,
 \item[.] for all objects $L$ and $M$ of $\cat$, 
we have $\chi(L\oplus M) = \chi(L)\chi(M)$, 
 \item[.] for all objects $L$ and $M$ of $\cat$ such that 
$\dim \ext^1_\cat(L,M) = 1$, we have 
$$
\chi(L)\chi(M) = \chi(B) + \chi(B'),
$$
\end{itemize}
where $B$ and $B'$ are the middle terms of `the' non-split 
triangles 
$$
L \fl B \fl M \fl \susp L \text{ and } 
M\fl B' \fl L \fl \susp M
$$
with end terms $L$ and $M$.
\end{defi}

Let $N$ be a finite-dimensional $B$-module and $e$ an element of $\ko$. 
We write $\gre(N)$ for the variety of submodules $N'$ of $N$
whose class in $\ko$ is $e$. 
It is a closed, hence projective, subvariety of the classical 
Grassmannian of subspaces of $N$. Let $\chi(\gre N)$ denote its 
Euler--Poincar\'e characteristic with respect to the \'etale 
cohomology with proper support. 
Let $\kosp$ denote the `split' Grothendieck group of $\modb$, i.e. 
the quotient of the free abelian group on the set of isomorphism classes $[N]$ of
finite-dimensional $B$-modules $N$, modulo the subgroup generated by all elements
$$
[N_1\oplus N_2] -[N_1] -[N_2] \text{.}
$$
We define a bilinear form 
$$
\langle \; ,\; \rangle : \kosp \times \kosp \gfl \mathbb{Z}
$$
by setting
$$
\langle N,N'\rangle  = [N,N']\!~-\!~^1[N,N']
$$
for all finite-dimensional $B$-modules $N$ and $N'$.
We define an antisymmetric bilinear form on $ \kosp$ by setting 
$$
\langle N,N'\rangle_a \;=\; \langle N,N'\rangle - \langle N',N\rangle 
$$
for all finite-dimensional $B$-modules $N$ and $N'$.
Let $T_1,\ldots, T_n$ be the pairwise non-isomorphic 
indecomposable direct summands of $T$ and, for 
$i= 1,\ldots,n$, let $S_i$ be the top of the projective 
$B$-module $P_i=FT_i$. The set $\{S_i, i = 1,\ldots,n\}$ 
is a set of representatives for the isoclasses of simple $B$-modules.

We need a lemma, the proof of which will be given in section~\ref{ss:welldef}.
\begin{lem}\label{lem:lf}
 For any $i=1,\ldots, n$, the linear form 
$\langle S_i,?\rangle_a \;:\, \kosp \fl \mathbb{Z}$ 
induces a well-defined form
$$
\langle S_i,?\rangle_a \;:\, \ko \fl \mathbb{Z}.
$$
\end{lem}
Let ind$\,\cat$ be a set of representatives for the isoclasses 
of indecomposable objects of $\cat$.
Define, as in \cite[6.1]{CK1}, a Caldero--Chapoton map,
$X^T_? : \text{ind}\,\cat~\fl~\mathbb{Q}(x_1,\ldots,x_n)$ by
$$
X^T_M = \left\{ \begin{array}{l} x_i \text{ if }M \simeq \susp T_i \\ 
\sum_e \chi(\gre FM)\prod_{i=1}^n x_i^{\langle S_i,e\rangle_a-\langle S_i,FM\rangle} \text{ else.} \end{array} \right.
$$
Extend it to a map $X^T_? : \cat \fl \mathbb{Q}(x_1,\ldots,x_n)$ 
by requiring that $X^T_{M\oplus N}~=~X^T_M X^T_N$. 
When there are no possible confusions, we often denote 
$X^T_M$ by $X_M$.
The main result of this article is the following
\begin{theo} \label{theo:mf}
The map $X^T_? : \cat \fl \mathbb{Q}(x_1,\ldots,x_n)$ 
is a cluster character.
\end{theo}

We will prove the theorem in section~\ref{ss:proof}, 
illustrate it by examples in section~\ref{ss:examples} 
and draw some consequences in section~\ref{ss:csq}.

\section{Index, coindex and Euler form}
\label{s:ind}

In the next two sections, our aim is to understand 
the exponents appearing in the definition of $X_M$. 
More precisely, for two objects 
$L$ and $M$ of $\cat$, we want to know how the exponents in $X_B$ 
depend on the choice of the middle term $B$ of a triangle 
with outer terms $L$ and $M$.

\subsection{Index and coindex}
\label{ss:ind}

\noindent Let $X$ be an object of $\cat$. 
Define its index\\ \noindent $\ind X\in\,\kop$ as follows. 
There exists a triangle (see [KR1])
$$
T^X_1 \fl T^X_0 \fl X \fl \susp T^X_1
$$
with $T^X_0$ and $T^X_1$ in add$\,T$. Define 
$\ind X$ to be the class $[FT^X_0] - [FT^X_1]$ 
in $\kop$.
Similarly, define the coindex of $X$,
denoted by $\coi X$, to be the class 
$[FT_X^0] - [FT_X^1]$ in $\kop$,
where 
$$
X \fl \susp^2 T_X^0 \fl \susp^2 T_X^1 \fl \susp X
$$
is a triangle in $\cat$ with $T_X^0,T_X^1\in\;$add$\,T$.\\

\begin{lem}\label{lem:ppties} We have the following properties:
\begin{enumerate}
\item The index and coindex are well defined.
\item $\ind X = - \coi\susp X$.
\item $\ind T_i = [P_i]$ and $\ind \susp T_i = -[P_i]$ where $P_i = FT_i$.
\item $\ind X - \coi X$ only depends on $FX\in$ mod$\,B$.
\end{enumerate}
\end{lem}

\begin{proof}
A right add$\,T$-approximation of an object $X$
of $\cat$ is a 
morphism $T'\stackrel{f}{\gfl} X$ with $T'\in\,$add$\,T$
such that any morphism $T''\gfl X$ with $T''\in\,$add$\,T$
factors through $f$. It is called minimal if, moreover,
any morphism $T' \stackrel{g}{\gfl}T'$ such that
$fg=f$ is an isomorphism. A minimal approximation
is unique up to isomorphism. \\ 
Assertions (2) and (3) are left to the reader.\\

\noindent (1) In any triangle of the form
$$
T^X_1 \fl T^X_0 \stackrel{f}{\fl} X \fl \susp T^X_1,
$$
the morphism $f$ is a right add$\,T$-approximation.
Therefore, any such triangle is obtained from one
where $f$ is minimal by adding a trivial triangle
$$
T' \fl T' \fl 0 \fl \susp T'
$$
with $T'\in\,$add$\,T$. The index is thus well-defined.
Dually, one can define left approximations and 
show that the coindex is well-defined.\\

\noindent (4) Let $T'$ be an object in add$\,T$. Take two triangles
$$
T^X_1 \fl T^X_0 \fl X \fl \susp T^X_1 \text{ and}
$$
$$
X \fl \susp^2 T_X^0 \fl \susp^2 T_X^1 \fl \susp X
$$
with $T^X_0$, $T^X_1$, $T_X^0$ and $T_X^1$ in add$\,T$.
Then, we have two triangles
$$
T^X_1 \oplus T' \fl T^X_0 \fl X\oplus \susp T' \fl \susp (T^X_1\oplus T') \text{ and}
$$
$$
X \oplus \susp T'\fl \susp^2 T_X^0 \fl \susp^2 (T_X^1\oplus T') \fl \susp X \oplus \susp^2 T'.
$$
We thus have the equality:
$$
\ind (X\oplus \susp T') - \coi (X\oplus \susp T') = \ind X - \coi X.
$$
\end{proof}

\begin{prop} \label{prop:ind}
Let $X\stackrel{f}{\fl} Z \stackrel{g}{\fl} Y \stackrel{\eps}{\fl} \susp X$ 
be a triangle in $\cat$. 
Take $C\in\cat$ (resp. $K\in\cat$) to be any lift of $\coke Fg$ (resp. $\ke Ff$). Then
$$
\begin{array}{rcl}
\ind Z & = & \ind X + \ind Y - \ind C - \ind \susm C \text{ and } \\
\coi Z & = & \coi X + \coi Y - \coi K - \coi \susp K.
\end{array}
$$
\end{prop}

\begin{proof} Let us begin with the equality for the indices. 
First, consider the case where $FC=0$. This means that the morphism
 $\eps$ belongs to the ideal $(\susp T)$. Take two triangles
$$
T_1^X \gfl T_0^X \gfl  X \gfl \susp T_1^X 
\text{ and } \;\; T_1^Y \gfl T_0^Y \gfl Y \gfl \susp T_1^Y
$$
in $\cat$, where the objects $T_0^X$, $T_1^X$, $T_0^Y$,  $T_1^Y$ 
belong to the subcategory add$\,T$. 
Since the morphism $\eps$ belongs to the ideal $(\susp T)$, 
the composition $T_0^Y \fl Y \stackrel{\eps}{\fl} \susp X$ 
vanishes. The morphism $T_0^Y \fl Y$ thus factors through $g$. 
This gives a commutative square
$$
\xymatrix{
T_0^X\oplus T_0^Y \dr \bas & T_0^Y  \bas \\
Z \dr & Y .
}
$$
Fit it into a nine-diagram
$$
\xymatrix{
T_1^X \dr \bas & Z' \dr \bas & T_1^Y \dr \bas & \susp T_1^X \\
T_0^X \dr \bas & T_0^X\oplus T_0^Y \dr \bas & T_0^Y \dr^0 \bas \bdr^0 \bgp & \susp T_0^X \\
X \dr \bas & Z \dr^g \bas & Y \dr^\eps \bas & \susp X \\
\susp T_1^X  & \susp Z' & \susp T_1^Y, & 
}
$$
whose rows and columns are triangles. 
Since the morphism $T_1^Y \fl  \susp T_1^X$ vanishes, 
the triangle in the first row splits, so that we have
$$
Z' \simeq T_1^X \oplus T_1^Y \text{ and }\; \ind Z = \ind X + \ind Y.
$$
Now, let us prove the formula in the general case. 
Let $FY \stackrel{a}{\gfl}M$ be a cokernel for $Fg$. 
Since the composition $F\eps\,Fg$ vanishes, 
the morphism $F\eps$ factors through $a$:
$$
\xymatrix{
FY \ddr ^{F\eps} \bdr_a & & F\susp X .\\
& M \hdrp_b & 
}
$$

Let $Y\stackrel{\alpha}{\gfl}C'$ be a lift of 
$a$ in $\cat$, and let $\beta$ be a lift of $b$.
The images under $F$ of the morphisms $\eps$ 
and $\beta\,\al$ coincide, therefore the morphism
$\beta\,\al - \eps$ belongs to the ideal $(\susp T)$.
Thus there exist an object $T'$ in add$\,T$ and two morphisms
$\al'$ and $\beta'$ such that the following diagram commutes:  
$$
\xymatrix{
Y \ddr ^\eps \bdr_{\left[^{\al}_{\al'} \right]} & & \susp X .\\
& C'\oplus \susp T' \hdr_{[\beta \; \beta']} & 
}
$$
Let $C$ be the direct sum $ C'\oplus \susp T'$. 

\noindent The octahedral axiom yields a commutative diagram

$$
\xymatrix{
& U \bas \ar@{=}[r] & U \bas & \\
X \dr \ar@{=}[d] & Z \dr \bas & Y \dr^\eps 
\bas^{^{\left[^{\al}_{\al'} \right]}} & \susp X \ar@{=}[d] \\
X \dr & V \dr^{\gamma'} \bas^\gamma & C \dr^{[\beta \; \beta']} \bas^{\gamma''} & \susp X \\
& \susp U \ar@{=}[r] & \susp U , & 
}
$$
whose two central rows and columns are triangles.
Due to the choice of $C$, the morphisms $\gamma'$, $\gamma''$, 
hence $\gamma$ belong to the ideal $(\susp T)$. 
We thus have the equalities:
\begin{eqnarray*}
\ind Y & = & \ind C + \ind U,  \\
\ind X & = & \ind V + \ind \susm C, \\
\ind Z & = & \ind V + \ind U, 
\end{eqnarray*}
giving the desired formula. Moreover, as seen in 
lemma~\ref{lem:ppties} (4), the sum
$\ind C +\ind \susm C =\ind C - \coi C$ 
does not depend on the particular choice of $C$.
Apply this formula to the triangle 
$$
\susm X \gfl \susm Z \gfl \susm Y \gfl X
$$
and use lemma~\ref{lem:ppties}(2) to obtain the formula 
for the coindices. Remark that the long exact sequence 
yields the equality of Coker$(\,-F\susm g)$ and Ker$\,Ff$.
\end{proof}

\subsection{Exponents.}\label{ss:exp}

\noindent We now compute the index and coindex in 
terms of the Euler form.

\begin{lem} \label{lem:ind}
Let $X\in\cat$ be indecomposable. Then
$$
\ind X = \left\{ \begin{array}{l} -[P_i] \text{ if } X \simeq \susp T_i
 \\ \\ \sum_{i=1}^n \langle FX,S_i \rangle [P_i] \text{ else,}  \end{array} \right.
$$
$$
\coi X = \left\{ \begin{array}{l} -[P_i] \text{ if } X \simeq \susp T_i 
\\ \\ \sum_{i=1}^n \langle S_i,FX\rangle [P_i] \text{ else.}  \end{array} \right.
$$
\end{lem}

\begin{proof}

Let $X$ be an indecomposable object in $\cat$, 
non-isomorphic to any of the $\susp T_i$'s. 
Take a triangle 
$$
T_1^X \stackrel{f}{\gfl} T_0^X \stackrel{g}{\gfl}  X \stackrel{\eps}{\gfl}\susp T_1^X 
$$ 
with the morphism $g$ being a minimal right add$\,T$-approximation, 
as defined in the proof of lemma~\ref{lem:ppties}. 
We thus get a minimal projective presentation
$$
P_1^X \gfl P_0^X \gfl FX \gfl 0
$$
where $P_i^X = FT_i^X, i=0,1$. 
For any $i$, the differential in the complex 
$$
0 \gfl (P_0^X, S_i) \gfl (P_1^X, S_i) \gfl \cdots
$$
vanishes. Therefore, we have
\begin{eqnarray*}
[FX,S_i] & = & [P_0^X, S_i] = [P_0^X : P_i],  \\
^1[FX,S_i] & = & [P_1^X, S_i] = [P_1^X : P_i],  \\
\langle FX,S_i \rangle & = & [\ind X : P_i]. 
\end{eqnarray*}
The proof for the coindex is analogous: 
We use a minimal injective copresentation of $FX$ 
induced by a triangle
$$
X \gfl \susp^2 T_X^0 \gfl \susp^2 T_X^1 \gfl \susp X.
$$
\end{proof}

\noindent Let us write $\underline{x}^e$ for 
$\prod_{i=1}^n x_i^{[e : P_i]}$ where $e \in\kop$ and 
$[e:P_i]$ is the $i$th coefficient of $e$ in the 
basis $[P_1],\ldots,[P_n]$. Then, by lemma~\ref{lem:ind}, 
for any indecomposable object $M$ in $\cat$, we have 
$$
X_M = \underline{x}^{-\coi M} \sum_e \chi(\gre FM) \prod_{i=1}^n x_i^{<S_i,e>_a}.
$$

\section{The antisymmetric bilinear form}
\label{s:abf}

 In this part, we give a positive answer to 
the first conjecture of \cite[6.1]{CK1} and
 prove that the exponents in $X_M$ are well defined. 
The first lemma is sufficient for this latter purpose, 
but is not very enlightening, whereas 
the second proof of theorem \ref{th:ko} gives us a 
better understanding of the antisymetric bilinear form. 
When the category $\cat$ is algebraic, this form is, 
in fact, the usual Euler form on the Grothendieck group of a 
triangulated category together with a t-structure whose 
heart is the abelian category mod$\,B$ itself. 

\subsection{The map $X^T$ is well defined}
\label{ss:welldef}

Let us first show that any short exact sequence in $\modb$ 
can be lifted to a triangle in $\cat$.

\begin{lem}\label{lem:ses}
Let $0\fl x\fl y \fl z \fl 0$ be a short exact sequence in $\modb$. 
Then there exists a triangle in $\cat$
$$
X \gfl Y \gfl Z \gfl \susp X
$$
whose image under $F$ is isomorphic to the given short exact sequence.
\end{lem}

\begin{proof}
 Let
$$
0 \gfl x \stackrel{i}{\gfl} y \stackrel{p}{\gfl} z \gfl 0
$$
be a short exact sequence in mod$\,B$. 
Let $X\stackrel{f}{\gfl} Y$ be a lift of the monomorphism 
$x \stackrel{i}{\gfl} y $ in $\cat$. 
Fix a triangle
$$
T_1^X \gfl T_0^X \gfl  X \gfl \susp T_1^X 
$$
and form a triangle
$$
X \gfl Y\oplus \susp T_1^X \gfl Z \stackrel{\eps}{\gfl} \susp X \text{ .}
$$
The commutative left square extends to a morphism of triangles
$$
\xymatrix{
X \ar@{=}[d] \dr & Y\oplus \susp T_1^X \bas^{[0\;1]} \dr & Z \basp \dr^\eps & \susp X \ar@{=}[d] \\
X \dr & \susp  T_1^X \dr & \susp  T_0^X \dr & \susp X.
}
$$
so that the morphism $\eps$ lies in the ideal $(\susp T)$. 
Therefore, the sequence
$$
0 \gfl x \stackrel{i}{\gfl} y \gfl FZ \gfl 0
$$
is exact, and the modules $FZ$ and $z$ are isomorphic.
\end{proof}

\noindent {\itshape Proof of lemma \ref{lem:lf}}.

Let $X$ be an object of the category $\cat$. 
Using section \ref{ss:exp} we have 
$$
\coi X - \ind X = \sum_{i=1}^n \langle S_i, FX\rangle_a \,[P_i]\text{ .} 
$$
Therefore, it is sufficient to show that the form
\begin{eqnarray}
\ko & \gfl  & \mathbb{Z} \nonumber \\
\text{[}FX\text{]} & \longmapsto & \coi X - \ind X \nonumber
\end{eqnarray}
is well defined. We already know that $\coi X - \ind X$ only depends on $FX$.
Take $0\fl x\fl y \fl z \fl 0$ to be a short exact sequence 
in $\modb$. Lift it, as in lemma~\ref{lem:ses}, to a triangle
$$
X \gfl Y \gfl Z \gfl \susp X \text{ in } \cat.
$$
 By proposition \ref{prop:ind}, we have
$$
\ind Y - \coi Y = (\ind X + \ind Z) - (\coi X + \coi Z) \nonumber
$$
which is the required equality. \qed

\begin{cor}
The map
$$
X^T_? \, : \cat \gfl \mathbb{Q}(x_1,\ldots,x_n)
$$
is well defined.
\end{cor}

\subsection{The antisymmetric bilinear form descends to the Grothendieck group.}

In this subsection, we prove a stronger result than in the previous one. 
This gives a positive answer to the first conjecture in \cite[6.1]{CK1}.

\begin{lem}\label{lem:dual}
 Let $T'$ be any cluster-tilting object in $\cat$.
We have bifunctorial isomorphisms
$$
\cat /_{(T')}(\susm X,Y) \simeq D(T')(\susm Y,X).
$$
\end{lem}

\begin{proof}
 Let $X$ and $Y$ be two objects of $\cat$, and let
$ T'_1 \gfl T'_0 \gfl X \stackrel{\eta}{\gfl} \susp T'_1$
be a triangle in $\cat$, with $T'_0$ and $T'_1$ in add$\,T'$.
Consider the morphism 
\begin{eqnarray*}
 \al : \; (T'_1,Y) & \gfl & (\susm X,Y) \\
 f & \longmapsto & f \circ \susm \eta.
\end{eqnarray*}
We have
$$
D(T')(\susm X,Y) \simeq D\im \al  \simeq \im D\al.
$$
Since the category $\cat$ is 2-Calabi--Yau, the dual 
of $\al$, $D\al$, is isomorphic to 
\begin{eqnarray*}
 \al' : \; (\susm Y,X) & \gfl & (\susm Y,\susp T'_1) \\
 g & \longmapsto & \eta \circ g.
\end{eqnarray*}
We thus have isomorphisms
\begin{eqnarray*}
D(T')(\susm X,Y) & \simeq & \im \al' \\
 & \simeq & (\susm Y,X)/\ke \al' \\
 & \simeq & \cat/_{(T')}(\susm Y,X).
\end{eqnarray*}

\end{proof}

\begin{theo}\label{th:ko}
The antisymmetric bilinear form $\langle \;,\;\rangle_a$ 
descends to the Grothendieck group $\ko$.
\end{theo}

\begin{proof}

Let $X$ and $Y$ be two objects in the category $\cat$. 
In order to compute $\langle FX, FY\rangle = [FX,FY] -\,\!^1[FX,FY]$,
let us construct a projective presentation in the following way. Let
$$
\susm X \stackrel{g}{\gfl} T_1^X \stackrel{f}{\gfl} T_0^X \gfl X
$$
be a triangle in $\cat$ with $T_0^X$ and $T_1^X$ being two objects in
the subcategory add$\,T$.
This triangle induces an exact sequence in $\modb$
$$
F\susm X \stackrel{Fg}{\gfl} FT_1^X \stackrel{Ff}{\gfl} FT_0^X \gfl FX \gfl 0,
$$
where $FT_0^X$ and  $FT_1^X$ are finite-dimensional projective $B$-modules.
Form the complex
$$
(\ast) \; \; \; \; \; 0 \gfl \homb ( FT_0^X, FY ) \gfl \homb (FT_1^X,FY) \gfl \homb(F\susm X, FY).
$$
Since the object $T$ is cluster-tilting in $\cat$, there are no morphisms 
from any object in add$\,T$ to any object in add$\,\susp T$. The complex
$(\ast)$ is thus isomorphic to the following one :
$$
0 \gfl \cat(T_0^X, Y ) \stackrel{f^\ast}{\gfl} \cat(T_1^X, Y ) \stackrel{g^\ast}{\gfl} \cat/_{\!(\susp T)}(\susm X,Y),
$$
where $f^\ast$ (resp. $g^\ast$) denotes the composition by $f$ (resp. $g$).
Therefore, we have 
\begin{eqnarray*}
\homb (FX,FY) & \simeq & \ke f^\ast \\
\ext^1_B(FX,FY) & \simeq & \ke g^\ast / \im f^\ast.
\end{eqnarray*}
We can now express the bilinear form as
\begin{eqnarray*}
\langle FX, FY\rangle & = & \dim \ke f^\ast - \dim \ke g^\ast + {\sf rk}\,f^\ast \\
 & = & [T_0^X,Y] - [T_1^X,Y] + {\sf rk}\,g^\ast , 
\end{eqnarray*}
with the image of the morphism $g^\ast$ being the quotient by the ideal $(\susp T)$
of the space of morphisms from $\susm X$ to $Y$, in $\cat$, 
which belong to the ideal $(T)$:
$$\im g^\ast = (T)/_{(\susp T)} (\susm X, Y ).$$
Similarily, using an injective copresentation given by 
a triangle of the form 
$$ X \gfl \susp ^2 T_X^0 \gfl \susp ^2 T_X^1 \stackrel{\beta}{\gfl}\susp X,$$
we obtain
$$
\langle FY, FX\rangle  = [Y,\susp^2 T_X^0] - [Y, \susp ^2 T_X^1 ] + {\sf rk}\,\beta_\ast,
$$
and $\im \beta_\ast = (\susp ^2 T)/_{(\susp T)} (Y,\susp X)$.
By lemma~\ref{lem:dual}, we have bifunctorial isomorphisms
$$
(T)/_{(\susp T)} (\susm X, Y ) \simeq D (\susp T)/_{(T)}(\susm Y,X) \simeq  D (\susp ^2 T)/_{(\susp T)} (Y,\susp X).
$$
Therefore, we have the equality
\begin{eqnarray*}
\langle FX, FY\rangle_a & = & [T_0^X,Y] - [T_1^X,Y]  - [Y,\susp^2 T_X^0] + [Y, \susp ^2 T_X^1 ] \\
 & = & [FT_0^X,FY] - [FT_1^X,FY] - [FY,F\susp^2 T_X^0] + [FY, F\susp ^2 T_X^1 ].
\end{eqnarray*}
Since $FT$ is projective and $F\susp^2 T$ in injective,
this formula shows that $\langle\;,\;\rangle_a$ 
descends to a bilinear form on the Grothendieck group $\ko$.
\end{proof}

\subsection{The antisymmetric bilinear form and the Euler form.}

In this subsection, assume moreover that the category $\cat$ 
is algebraic, as in \cite[section 4]{KR1}: There exists a $k$-linear
 Frobenius category with split idempotents $\ec$ 
whose stable category is $\cat$. 
Denote by $\mc$ the preimage, in $\ec$, of add$\,T$ 
via the canonical projection functor. 
The category $\mc$ thus 
contains the full subcategory $\pc$ of $\ec$ 
whose objects are the projective objects in $\ec$, 
and we have $\ms =\;$add$\,T$.
Let $\modm$ be the category of $\mc$-modules, 
i.e. of $k$-linear contravariant functors from $\mc$ to the category of 
$k$-vector spaces. 
The category $\modms$ of finitely presented 
$\ms$-modules is identified with 
the full subcategory of $\modm$ 
of finitely presented $\mc$-modules vanishing on $\pc$. 
This last category is equivalent to the abelian category $\modb$ 
of finitely generated $B$-modules. 
Recall that the perfect derived category $\perm$ is the 
full triangulated subcategory of the derived
category of ${\sf \mathcal{D}}\modm$ generated by 
the finitely generated projective $\mc$-modules.
Define $\perms$ to be the full subcategory of $\perm$ 
whose objects $X$ satisfy the following conditions:
\begin{itemize}
\item[.] for each integer $n$, the finitely presented $\mc$-module ${\sf H}^n X$ belongs to $\modms$,
\item[.] the module ${\sf H}^n X$ vanishes for all but finitely many $n \in \z$.
\end{itemize}
It can easily be shown that $\perms$ is a 
triangulated subcategory of $\perm$. 
Moreover, as shown in \cite{Tab}, the canonical t-structure on 
${\sf \mathcal{D}}\modm$ induces a t-structure on $\perms$, 
whose heart is the abelian category $\modms$.

\noindent The following lemma shows that the Euler form 
\begin{eqnarray*}
{\sf K}_0 \left(\perms\right) \times {\sf K}_0 \left(\perms\right)  & \gfl & \z \\
 ([X],[Y]) & \longmapsto & \langle[X],[Y]\rangle = \sum_{i\in \z} (-1)^i \dim \perms\left(X, \susp^i Y\right)
\end{eqnarray*}
is well defined.

\begin{lem}\label{lem:euler}
Let $X$ and $Y$ belong to $\perms$.
Then the vector spaces 
$\perms\left(X, \susp^i Y\right)$ are finite dimensional 
and only finitely many of them are non-zero.
\end{lem}

\begin{proof}
Since $X$ belongs to $\perm$, we may assume 
that it is representable: There exists 
$M$ in $\mc$ such that $X = M\hat\,$. 
Moreover, the module ${\sf H}^n Y$ 
vanishes for all but finitely many 
$n \in \z$. We thus may assume $Y$ to 
be concentrated in degree $0$. 
Therefore, the space 
$\perms\left(X, \susp^i Y\right) = 
\perms(M\hat\, ,\susp^i {\sf H}^0 Y)$ 
vanishes for all non-zero $i$. For 
$i=0$, it equals 
\begin{eqnarray*}
\homph_{\mc}\left(M\hat\, ,{\sf H}^0 Y
\right) & = & {\sf H}^0 Y(M) \\
 & = & \homph_{\ms}\left( 
\ms ( ?,M), {\sf H}^0 Y \right).
\end{eqnarray*}
this last space being finite dimensional.

\end{proof}

This enables us to give another proof of theorem \ref{th:ko}. 
This proof is less general than the previous one, but is nevertheless
much more enlightening.\\

\noindent {\itshape Proof of theorem \ref{th:ko}}.
Let $X$ and $Y$ be two finitely presented $\ms$-modules, 
lying in the heart of the t-structure on $\perms$. 
We have: 
\begin{eqnarray}
\langle[X],[Y]\rangle & = & \sum_{i\in \z} (-1)^i \dim \perms\left(X, \susp^i Y\right) \nonumber \\
 & = & \sum_{i = 0}^3 (-1)^i \dim \perms\left(X, \susp^i Y\right) \\
 & = & \dim\perms(X,Y) - \dim\perms(X,\susp Y) \nonumber \\
 & & \;\;\;\;\;\; \;\;\;\;\;\; 
+ \dim\perms(X,\susp^2 Y) -\dim\perms(X,\susp^3 Y) \nonumber \\
 & = &  \dim\perms(X,Y) - \dim\perms(X,\susp Y)  \\
& &  \;\;\;\;\;\; \;\;\;\;\;\;
+ \dim\perms(Y,X) -\dim\perms(Y,\susp X)  \nonumber\\
 & = & \dim\homph_{\ms}(X,Y) - \dim\ext^1_{\ms}(X,Y) \nonumber \\
 & & \;\;\;\;\;\; \;\;\;\;\;\; +\dim\homph_{\ms}(Y,X) -\dim\ext^1_{\ms}(Y,X) \nonumber \\
 & = & \langle [X],[Y]\rangle_a \nonumber
\end{eqnarray}
where the classes are now taken in $\ko$. 
Equalities $(1)$ and $(2)$ are consequences of the 
3-Calabi--Yau property of the category $\perms$, 
cf.~\cite{KR1}.
\qed

\section{Dichotomy}
\label{s:dicho}
Our aim in this part is to study the coefficients appearing in the 
definition of $X_M$. In particular, we will prove that the phenomenon 
of dichotomy proved in \cite{CK2} (see also \cite{GLS3}) remains true
 in this more general setting.

Recall that we write $\underline{x}^e$ for 
$\prod_{i=1}^n x_i^{[e : P_i]}$ where $e \in\kop$ and 
$[e:P_i]$ is the $i$th coefficient of $e$ in the 
basis $[P_1],\ldots,[P_n]$.

\begin{lem}\label{lem:dec}
For any $M\in\cat$, we have
$$
X_M =  \underline{x}^{-\coi M} \sum_e \chi(\gre FM) \prod_{i=1}^n x_i^{\langle S_i,e\rangle_a}.
$$
\end{lem}

\begin{proof}
We already know that this formula holds for indecomposable 
objects of $\cat$, cf. section~\ref{ss:exp}. Let us 
prove that it still holds for decomposable objects,
by recursion on the number of indecomposable 
direct summands.

Let $M$ and $N$ be two objects in $\cat$. As shown in \cite{CC}, we have 
$$
\chi\left(\grg F(M\oplus N) \right) = \sum_{e+f=g} \chi\left(\gre FM\right) \chi\left( \grf FN \right).
$$
Therefore, we have $X_{M\oplus N} = X_M X_N =$
$$
\left( \underline{x}^{-\coi M} \sum_e \chi(\gre FM) \prod_{i=1}^n x_i^{<S_i,e>_a}\right) \left(
\underline{x}^{-\coi N} \sum_f \chi(\gre FN) \prod_{i=1}^n x_i^{<S_i,f>_a} \right)
$$
$$
=  \underline{x}^{-(\coi M + \coi N)} \sum_g \sum_{e+f=g}  \chi\left(\gre FM\right) \chi\left( \grf FN \right)\prod_{i=1}^n x_i^{<S_i,e+f>_a}
$$
$$
=  \underline{x}^{-\coi (M\oplus N)} \sum_g\chi\left(\grg F(M\oplus N) \right)\prod_{i=1}^n x_i^{<S_i,g>_a}
$$
\end{proof}

\begin{lem}
Let $M\stackrel{i}{\gfl}B\stackrel{p}{\gfl}L\stackrel{\eps}{\gfl}\susp M $ 
be a triangle in $\cat$, and let $U\stackrel{i_U}{\gfl} M$ and
$V\stackrel{i_V}{\gfl} L$ be two morphisms whose images 
under $F$ are monomorphisms. Then the following 
conditions are equivalent:
\begin{itemize}
\item[i)] There exists a submodule $E\subset FB$ such that \\ $FV = (Fp)E$ and $FU=(Fi)^{-1}E$,
\item[ii)] There exist two morphisms $e : \susm V \gfl U$ and $f:\susm L \gfl U$ such that
	\begin{itemize}
	\item[a)] $(\susm\eps)(\susm i_V)  =  i_U e$
	\item[b)] $e \in (T)$  
	\item[c)] $i_U f  - \susm\eps \; \in (\susp T)$.
	\end{itemize}
\item[iii)] Condition ii) where, moreover, $e = f\susm i_V $.
\end{itemize}
\end{lem}

\noindent The following diagrams will help the reader parse the conditions:
$$
\xymatrix{
F\susm L \dr^{F\susm \eps} & FM \dr^{Fi} & FB \dr^{Fp} & FL & \\ 
& FU \ham \dr & E \ham \dr & FV \ham \dr & 0,
}
$$
$$
\xymatrix{
\susm L \dr^{\susm \eps} \bdr^f & M \\
\susm V \ha^{\susm i_V} \dr_e & U \ha_{i_U}.
}
$$

\begin{proof}
Assume condition ii) holds. Then, by a), there exists a morphism of triangles
$$
\xymatrix{
\susm L \dr^{\susm\eps} & M \dr^i & B \dr^p & L \\
\susm V \ha_{\susm i_V} \dr ^e & U \ha_{i_U} \drp & W \drp \hap_j & V \ha_{i_V}
}
$$
Take $E$ to be the image of the morphism $Fj$. The morphism $e$ factors through add$\,T$, 
so that we have $F\susp e=0$ and the functor $F$ induces a 
commutative diagram 
$$
\xymatrix{
 F\susm L \dr^{F\susm\eps} \ar@{->}[ddr]_{Ff} &   FM \dr^{Fi} & FB \dr^{Fp} & FL \dr^{F\eps} & F\susp M \\
  & & E \ham \ar@{-->>}[dr] & & \\
 &   FU \hham_{Fi_U} \dr \hdr & FW \dr \hae & FV \hham_{Fi_V} \dr & 0
}
$$ 
whose rows are exact sequences. 
It remains to show that $FU= (Fi)^{-1}E$. 

We have $FU \subset (Fi)^{-1}E$ since $(Fi)(Fi_U)$ factors 
through the monomorphism $E\fl FB$.  
The existence of the morphism $Ff$ shows, 
via diagram chasing, the converse inclusion.\\

Conversely, let $E\subset FB$ be such that $FV = (Fp)E$ 
and $FU=(Fi)^{-1}E$. In particular, $FU$ contains 
$\ke Fi = \im F\susm\eps$ so that 
 $F\susm\eps$ factors through $Fi_U$. 
This gives us the morphism $f$, satisfying condition c). 
Define the morphism $e$ as follows. There exists a triangle
$$
T_1 \gfl T_0 \gfl V \gfl \susp T_1,
$$
where $T_1,T_0$ belong to add$\,T$. 
Applying the functor $F$ to this triangle, we get an 
epimorphism $FT_0 \fl FV$ with  $FT_0$ projective.
This epimorphism thus factors through the surjection 
$E\fl FV$, and composing it with $E\fl FB$ gives a commutative square
$$
\xymatrix{
FT_0 \dr \bas & FV \bas \\
FB \dr & FL.
}
$$
Since $\cat(T,\susp T)=0$, this commutative square lifts
 to a morphism of triangles
$$
\xymatrix{
\susm V \bas \dr & T_1 \bas \dr & T_0 \bas \dr & V \bas \\
\susm L \dr & M \dr & B \dr & L.
}
$$ 
The morphism $T_1 \fl M$ thus induced, factors through
 the morphism $U\fl M$. Indeed, we have $FU=(Fi)^{-1}E$ and the 
following diagram commutes :
$$
\xymatrix{
& FM \ddr & & FB \\ 
FT_1 \hdr \ar@{-}[r] & \dr & FT_0 \hdr \bdr & \\
& FU \ddr \ar@{^{(}->}[uu] & & E .\ar@{^{(}->}[uu] 
}
$$
The morphism $e$ is then given by the composition $\susm V \gfl T_1 \gfl U$.\\

\noindent Let us show that condition ii) implies condition iii).
By hypothesis, we have
$$i_Ue = (\susm \eps) (\susm i_V) $$
and
$$ i_Uf\susm i_V \equiv (\susm \eps) (\susm i_V) \bmod (\susp T). $$
Therefore, the morphism $i_U\left( f\susm i_V - e \right)$ 
belongs to the ideal $(\susp T)$. 
The morphism $Fi_U$ is a monomorphism, 
so that the morphism $h := f\susm i_V - e$ 
lies in $(\susp T)$.
There exists a morphism $\susm L \stackrel{l}{\gfl}U$ 
such that $ h = l \susm i_V $ :
$$
\xymatrix{
\susm C \dr^{\in ( T)} \bdr_0 & \susm V \dr^{\susm i_V} 
\bas_h^{\in (\susp T)} & \susm L \dr^c \ar@{-->}[dl]^l & C \\
 & U & &.
}
$$
Since the morphism $\susm C \fl \susm V$ 
lies in the ideal $(T)$, 
there exists a morphism of triangles
$$
\xymatrix{
\susm C \dr \bas & \susm V \ar@{=}[d] \dr & \susm L \dr^c \bas_v & C \bas\\
T_V^1 \dr^u & \susm V \dr & \susp T_V^0 \dr & \susp T_V^1 . 
}
$$
The composition $l\susm i_V$ belongs to the ideal $(\susp T)$, 
so that the composition $l(\susm i_V)u$ vanishes. 
We thus have a morphism of triangles
$$
\xymatrix{
T_V^1 \dr^u \bas & \susm V \dr \bas_{\susm i_V} & \susp T_V^0 \dr \bas_w & \susp T_V^1 \bas \\ 
\susm C' \dr & \susm L \dr^l & U \dr & C'.
}
$$
Therefore, we have $(\susm i_V)(l - wv) = 0$, and there exists
a morphism $C\stackrel{l'}{\gfl} U$ such that
$l - wv = l'c$. The morphism $l_0 = l - l'c $ 
thus factors through $\susp T_1$.
Put $f_0 = f - l_0$. We have
$$ f_0 \susm i_V = f \susm i_V - l \susm i_V + l'c\susm i_V = e $$
and
\begin{eqnarray*}
i_U f_0 & = & i_Uf -i_U l_0 \\
 & \equiv & i_U f \mod (\susp T) \\
 & \equiv & \susm \eps \mod (\susp T).
\end{eqnarray*}

\end{proof}

\begin{prop}\label{prop:dicho}
Let $L,M\in\cat$ be such that $\dim \cat(L,\susp M) = 1$. 
Let 
$$
\Delta : M \stackrel{i}{\gfl} B \stackrel{p}{\gfl} L \stackrel{\eps}{\gfl}\susp M
$$
$$
 \!\!\!\!\!\!\!\!\!\!\text{and } \Delta' : L \stackrel{i'}{\gfl} B' \stackrel{p'}{\gfl} M 
\stackrel{\eps'}{\gfl}\susp L
$$
be non-split triangles.
Then conditions i) to iii) hold for the triangle $\Delta$ 
if and only if they do not for the triangle $\Delta'$.
\end{prop}

\begin{proof}

Define maps
\begin{eqnarray*}
(\susm L, U) \oplus (\susm L, M) & \stackrel{\al}{\gfl} & \cat/(T)\left( \susm V, U \right) \oplus (\susm V, M ) \oplus  \cat/(\susp T)\left(\susm L, M \right)
  \\
(f, \eta ) & \longmapsto & (f\susm i_V, i_U f \susm i_V - \eta \susm i_V , i_U f - \eta )  \\
\text{ and} \;\;\;\;\;\;\;\;\;\;\;\;\;\;\;\;\;\;\;\;\;\;\;\;\;\;\;\;\;\;\;\;\;\;\;\;\;\;\;\;\;\;\;\;\;
\;\;\;\;\;\;\;\;\;\;\; \;\;\;& & \\
(\susm U, L ) \oplus (\susm M, L ) & \stackrel{\al'}{\longleftarrow} & (T)(\susm U, V) \oplus (\susm M, V) \oplus (\susp T)(\susm M, L) \\
(i_V e' + g' \susm i_U + i_V f' \susm i_U, -g' -i_V f') & \longmapsfrom & (e',f',g'). 
\end{eqnarray*}

Since the morphism space $\cat(L,\susp M)$ is one-dimensional, 
the morphism $\eps$ satisfies condition iii) if and only if the composition
$$
\xymatrix{
\beta : \ke \al \;\; \ar@{^{(}->}[r] & (\susm L, U) \oplus (\susm L, M) \ar@{->>}[r] & (\susm L, M)
}
$$
does not vanish.
Assume condition iii) to be false for the triangle $\Delta$. 
This happens if and only if the morphism $\beta$ vanishes, 
if and only if its dual $D\beta$ vanishes. 
Since the category $\cat$ is 2-Calabi--Yau, 
lemma~\ref{lem:dual} implies that
the morphism $D\beta$ is isomorphic to the morphism:

$$
\xymatrix{
\beta' : (\susm M,L)\;\; \ar@{^{(}->}[r] & (\susm U,L)\oplus(\susm M,L) \ar@{->>}[r] & \coke \al'.
}
$$
Therefore, $\beta'(\susm\eps) = 0$
is equivalent to $\susm \eps$ being in $\im\al'$, which is 
equivalent to the existence of three mophisms $e',f',g'$ 
as in the diagram
$$
\xymatrix{
\susm M \dr^{g'} \bdr^{f'} & L \\
\susm U \ha^{\susm i_U} \dr_{e'} & V \ha_{i_V}
}
$$
such that
$$
\left\{ \begin{array}{l} 
e'\in(T) \\
g' \in (\susp T) \\
\susm \eps' = i_Vf' + g' \\
i_V e' = (\susm \eps')(\susm i_U).
\end{array} \right.
$$
We have thus shown that condition iii) does not hold for the
triangle $\Delta$
if and only if condition ii) holds for the triangle $\Delta'$.
\end{proof}

\section{The multiplication formula}
\label{s:mf}

We use sections \ref{s:ind} and \ref{s:dicho} to prove the multiplication formula, 
and apply it to prove conjecture 2 in \cite{CK1}.

\subsection{Proof of theorem \ref{theo:mf}}\label{ss:proof}

We use the same notations as in the statement of theorem \ref{theo:mf}.

\noindent Define, for any classes $e,f,g$ in the Grothendieck group $\ko$, 
the following varieties 
$$
\begin{array}{lcl}
 X_{e,f} & = & \{ E \subset FB \text{ s.t. } [(Fi)^{-1}E] = e \text{ and } [(Fp)E]  = f \} \\
 Y_{e,f} & = & \{ E \subset FB' \text{ s.t. } [(Fi')^{-1}E] = f \text{ and } [(Fp')E]  = e \} \\
 X_{e,f}^g & = & X_{e,f} \cap \grg (FB) \\
 Y_{e,f}^g & = & Y_{e,f} \cap \grg (FB').
\end{array}
$$
We thus have 
$$
\grg(FB) = \coprod_{e,f}X_{e,f}^g \text{ and } \grg(FB') = \coprod_{e,f}Y_{e,f}^g.
$$
Moreover, we have
\begin{eqnarray*}
\chi\left( \gre(FM)\times\grf(FL)\right) & = &
\chi\left( X_{e,f} \sqcup Y_{e,f}\right) \\
& = & \chi\left( X_{e,f} \right) + \chi\left( Y_{e,f} \right) \\
& = & \sum_g \left( \chi\left( X_{e,f}^g \right) + \chi\left(Y_{e,f}^g\right)\right).
\end{eqnarray*}
where the first equality is a consequence of 
the dichotomy phenomenon as follows:
Consider the map 
$$
X_{e,f} \sqcup Y_{e,f} \gfl \gre(FM)\times\grf(FL) 
$$
which sends a submodule $E$ of $FB$ to the pair of 
submodules $\left( (Fi)^{-1}E , (Fp)E \right)$. 
By proposition~\ref{prop:dicho}, it is surjective, 
and, as shown in~\cite{CC}, its fibers are affine 
spaces.

\begin{lem}\label{lem:exp}
 Let $e,f$ and $g$ be classes in 
$\text{\rm K}_0(\operatorname{mod}\operatorname{End}_{\cat}(T))$. 
Assume that $X_{e,f}^g$ is non-empty. Then, we have
$$
\sum\langle S_i,g\rangle_a [P_i] - \coi B = 
\sum\langle S_i, e+f\rangle_a -\coi M -\coi L.
$$
\end{lem}

\begin{proof}
 Let $E$ be a submodule of $FB$ in $X_{e,f}^g$. 
Let $U\stackrel{i_U}{\gfl}M$ and $V \stackrel{i_V}{\gfl}L$ 
be two morphisms in the category $\cat$ such that
$FU\simeq (Fi)^{-1}E$, $FV\simeq (Fp)E$ and 
the images of $i_U$ and $i_V$ in $\modb$
are isomorphic to the inclusions of 
$FU$ in $FM$ and $FV$ in $FL$ respectively.
Let $K\in\cat$ be a lift of the kernel of $Fi$.
By proposition~\ref{prop:ind}, the following 
equality holds: 
$$
(1) \;\;\;\; \coi B = \coi M + \coi L - \coi K -\coi (\susp K).
$$
By diagram chasing, the kernel of $Fi$ is also
a kernel of the induced morphism from
$FU$ to $E$. Therefore, in $\ko$, we have
$$
(2) \;\;\;\; g = e + f - [FK].
$$
We have the following equalities: 
\begin{eqnarray*}
\sum \langle S_i, FK \rangle_a [P_i]  & = & 
\coi K - \ind K \text{ (by lemma~\ref{lem:ind})} \\
 & = & \coi K + \coi (\susp K) \text{ (by lemma~\ref{lem:ppties}).}
\end{eqnarray*}
Equality $(2)$ thus yields
$$
(3) \;\;\;\; \sum \langle S_i, g \rangle_a [P_i] =
\sum \langle S_i, e+f \rangle_a [P_i] - \coi K -\coi (\susp K).
$$
It only remains to sum equalities $(1)$ and $(3)$ to finish the proof.
\end{proof}

\noindent {\itshape Proof of theorem \ref{theo:mf}.}

Using lemma~\ref{lem:dec}, we have 
\begin{eqnarray*}
X_M X_L & = & \underline{x}^{-\coi M-\coi L}\sum_{e,f}\chi(\gre FM)\chi(\grf FL)
\prod_{i=1}^n x_i^{\langle S_i,e+f\rangle_a},  \\
X_B & = & \underline{x}^{-\coi B}\sum_g \chi(\grg FB)
\prod_{i=1}^n x_i^{\langle S_i,g\rangle_a} \text{ and} \\
X_{B'} & = & \underline{x}^{-\coi B'}\sum_{g} \chi(\grg FB')
\prod_{i=1}^n x_i^{\langle S_i,g\rangle_a}. 
\end{eqnarray*}
Therefore
\begin{eqnarray*}
 X_M X_L & = & \underline{x}^{-\coi M -\coi L}\sum_{e,f}
\chi\left( \gre (FM) \right) \chi \left(\grf (FL) \right)
\prod x_i^{\langle S_i, e+f \rangle_a} \\
 & = & \underline{x}^{-\coi M -\coi L}\sum_{e,f,g}
\left( \chi \left( X_{e,f}^g \right) + \chi \left( Y_{e,f}^g\right) \right)
\prod x_i^{\langle S_i, e+f \rangle_a} \\
 & = &
\underline{x}^{-\coi B} \sum_{e,f,g} \chi \left( X_{e,f}^g \right)
\prod x_i^{\langle S_i, g \rangle_a} \\
 & & \;\;\;\;\;\;\;\;\;\;\;\;\;\;\;\;\;\;\;\; + 
\underline{x}^{-\coi B'} \sum_{e,f,g} \chi \left( Y_{e,f}^g \right)
\prod x_i^{\langle S_i, g \rangle_a} \\
 & = & \underline{x}^{-\coi B} \sum_{g} 
\chi\left( \grg(FB) \right) \prod x_i^{\langle S_i, g \rangle_a} \\
 & & \;\;\;\;\;\;\;\;\;\;\;\;\;\;\;\;\;\;\;\; + \underline{x}^{-\coi B'} \sum_{g} 
\chi\left( \grg(FB') \right) \prod x_i^{\langle S_i, g \rangle_a} \\
 & = & X_B + X_{B'}.
\end{eqnarray*}
\qed

\subsection{Consequences}\label{ss:csq}

Let $Q$ be a finite acyclic connected quiver, 
and let $\cat$ be the cluster category 
associated to $Q$.

An object of $\cat$ without self-extensions 
is called rigid.
An object of $\cat$ is called basic if 
its indecomposable direct summands are 
pairwise non-isomorphic.
For a basic cluster-tilting object $T$ 
of $\cat$, let $Q_T$ denote the quiver 
of End$\,(T)$, and $\aqt$ the associated 
cluster algebra.

\begin{prop}\label{prop:csq}
A cluster character $\chi$ on $\cat$ 
with values in $\mathbb{Q}(x_1,\ldots,x_n)$ which sends 
a basic cluster-tilting object $T$ of 
$\cat$ to a cluster of $\aqt$ sends any 
cluster-tilting object $T'$ of $\cat$ 
to a cluster of $\aqt$, and 
any rigid indecomposable object to a 
cluster variable.
\end{prop}

\begin{proof}
Since the tilting graph of $\cat$ is 
connected, cf.~\cite[proposition 3.5]{BMRRT},
we can prove the first part of the 
proposition by recursion on the minimal 
number of mutations linking $T'$ to $T$. 
Let $T''= T''_1 \oplus \cdots \oplus T''_n$ 
be a basic cluster-tilting object, 
whose image under $\chi$ is a cluster of $\aqt$.
Assume that $T'= T'_1\oplus T''_2 \oplus 
\cdots \oplus T''_n$ is the mutation in direction 
1 of $T''$. Since $\chi$ is a cluster character, 
it satisfies the multiplication formula, and 
theorem 6.1 of~\cite{BMR2} shows that 
the mutation, in direction~1, of the cluster 
$\left( \chi(T''_1), \ldots, \chi(T''_n) \right)$ 
is the cluster $\left( \chi(T'_1), 
\chi(T''_2), \ldots, \chi(T''_n) \right)$. 
We have thus proved that the image under 
$\chi$ of any cluster-tilting object is a 
cluster. 
It is proved in~\cite[proposition 3.2]{BMRRT} 
that any rigid indecomposable object 
of $\cat$ is a direct summand of a 
basic cluster-tilting object. 
Therefore, the image under $\chi$ of 
any rigid indecomposable object is 
a cluster variable of $\aqt$. 
\end{proof}

\noindent {\itshape Remark}: As a corollary of 
the proof of proposition~\ref{prop:csq}, 
a cluster character is characterised, 
on a set of representatives for the isoclasses 
of indecomposable rigid objects of $\cat$ by 
the image of each direct summand of any 
given cluster-tilting object. 
In fact, using~\cite[1.10]{BIRS}, 
this remains true in the more general context 
of~\cite{BIRS}: 
Let $\cat$ be a Hom-finite triangulated 2-Calabi--Yau 
category having maximal rigid objects without loops nor 
strong 2-cycles. Denote by $n$ the number of non-isomorphic 
indecomposable direct summands of any maximal rigid object.

\begin{lem}\label{Schroer}
Let $\chi_1$ and $\chi_2$ be two cluster characters on $\cat$ 
with values in $\mathbb{Q}(x_1,\ldots,x_n)$. Assume that $\chi_1$ 
and $\chi_2$ coincide on all indecomposable direct summands 
of a cluster-tilting object $T$ in $\cat$.
Then $\chi_1$ and $\chi_2$ coincide on all 
direct summands of the cluster-tilting objects in $\cat$ 
which are obtained from $T$ by a finite sequence of mutations.
\end{lem}

The following corollary was conjectured for 
the finite case in~\cite{CK1}: Let $\cat$ be the cluster 
category of the finite acyclic quiver $Q$.
\begin{cor}\label{cor:csq}
Let $T$ be any basic 
cluster-tilting object in $\cat$, and let 
$Q_T$ denote the quiver of End$\,(T)$. 
Denote by $\tc$ a 
set of representatives for the isoclasses 
of indecomposable rigid objects of $\cat$.
Then $X^T$ induces a bijection 
from the set $\tc$ 
to the set of cluster variables of the 
associated cluster algebra $\aqt$, sending 
basic cluster-tilting objects to clusters.
\end{cor}

\begin{proof}
In view of theorem~\ref{theo:mf}, 
proposition~\ref{prop:csq} shows that
the map $X^T$ sends rigid indecomposable 
objects to cluster variables and 
cluster-tilting objects to clusters. 
It remains to show that it induces 
a bijection. This follows 
from~\cite[theorem 4]{CK2}, 
where it is proved 
for the Caldero-Chapoton map $X^{kQ}$.

As in the proof of proposition~\ref{prop:csq}, 
we proceed by induction on the minimal 
number of mutations linking $T$ to $kQ$.

Let $T'$ be a basic cluster-tilting object 
such that the map $X^{T'}$ induces a bijection 
from the set $\tc$ to the set of cluster 
variables. Assume that $T$ is the mutation 
in direction~1 of $T'$. 
Denote by $f$ the canonical isomorphism from 
$\aqtp$ to $\aqt$. Theorem~6.1 of~\cite{BMR2} 
shows that the two cluster characters 
$X^T$ and $f\circ X^{T'}$ coincide on the 
indecomposable direct summands of $\susp T$. 
Therefore, they coincide on all rigid objects and
the map $X^T$ also induces a bijection. 
\end{proof}

\noindent {\itshape Remark}: 
We have shown that, for any basic 
cluster-tilting object $T$, we have a 
commutative diagram
$$
\xymatrix{
 & \tc \ar@{->}[dl]\ar@{->}[dr]^{X^{T}} & \\
\aq & & \aqt \ar@{->}[ll]_{\simeq}
}
$$
where the arrow on the left side is the 
Caldero--Chapoton map.

\section{Examples}\label{ss:examples}

\subsection{The cluster category $\cat_{A_4}$}

The Auslander--Reiten quiver of $\cat_{A_4}$ is
$$
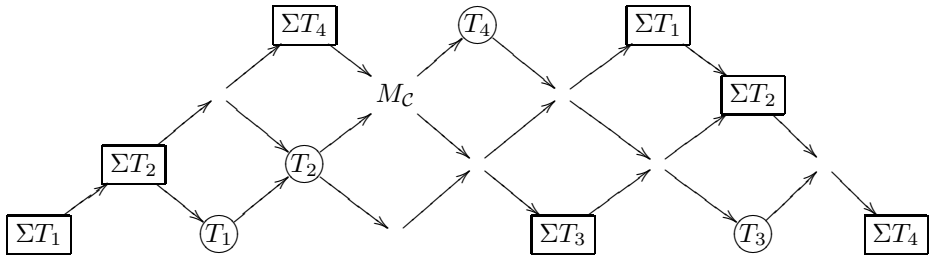
\xymatrix@-1pc{
 & & & *+[F]{\susp T_4}  
\bdr & & *+[o][F]{T_4} \bdr & &*+[F]{\susp T_1} \bdr & & & \\
 & &  \bdr \hdr & & M_{\cat} \bdr \hdr & &  \bdr \hdr & & *+[F]{\susp T_2}\bdr & & \\
 & *+[F]{\susp T_2} \bdr \hdr & & *+[o][F]{T_2} \bdr \hdr & &  \bdr \hdr & &  \bdr \hdr & &  \bdr & \\
 *+[F]{\susp T_1}\hdr & & *+[o][F]{T_1} \hdr & &  \hdr & & *+[F]{\susp T_3}\hdr & 
& *+[o][F]{T_3} \hdr & & *+[F]{\susp T_4}.
}
$$
The object $T:=T_1\oplus T_2 \oplus T_3 \oplus T_4$ is cluster-tilting.
Indeed, it is obtained from the image of the $kQ$-projective module $kQ$
in $\cat_{A_4}$ by the mutation of the third vertex. 

The quiver of $B=\,$End$_{\cat_{A_4}}(T)$ is 
$$
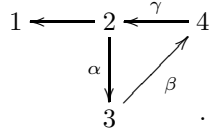
\xymatrix{
1 & 2 \ar@{->}[l] \bas_\al & 4 \ar@{->}[l]_\gamma \\
& 3 \hdr_\beta & .
}
$$
with relations $\beta\al = \gamma\beta = \al\gamma = 0$.
For $i=1,\ldots,n$, let $P_i$ be the image of $T_i$ in $\modb$, 
let $I_i$ be the image of $\susp^2 T_i$ and let $S_i$ be the simple 
top of $P_i$. Let $M$ be the finite-dimensional $B$-module given by: 

$$
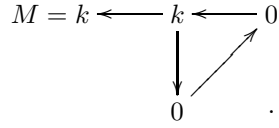
\xymatrix{
M = k & k \ar@{->}[l] \bas & 0 \ar@{->}[l]  \\
 & 0 \hdr & .
}
$$

The shape and the relations of the AR--quiver of $B$ 
are obtained from the ones of $\cat_{A_4}$
by deleting the vertices corresponding to the objects $\susp T_i$
and all arrows ending to or starting from these vertices. 
$$
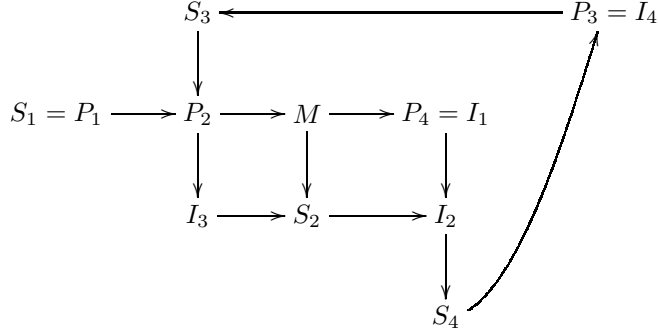
\xymatrix{
 & S_3 \bas & & & P_3 = I_4 \ar@{->}[lll] \\
 S_1 = P_1 \dr & P_2 \bas \dr & M \bas \dr & P_4 = I_1 \bas & \\
 & I_3 \dr & S_2 \dr & I_2 \bas & \\
 & & & S_4 \ar@(r,u)@{->}[ruuu] &
}
$$
Let $M_{\cat}$ be an indecomposable lift of $M$ in $\cat_{A_4}$.
The triangles
$$
T_3 \gfl T_2 \gfl M_{\cat} \gfl \susp T_3
\;\; \text{ and } \;\;
T_1 \gfl T_4 \gfl \susm M_{\cat} \gfl \susp T_1
$$
allows us to compute the index and coindex of $M_{\cat}$: 
\begin{eqnarray*}
 \ind M_{\cat} & = & [P_2] - [P_3] \\
\coi M_{\cat} & = & [P_1] - [P_4].
\end{eqnarray*}
Up to isomorphism, the submodules of $M$ are 
$0$, the simple $S_1$, and $M$ itself.
We thus have 
$$
X_{M_{\cat}} = \frac{x_4 x_2 + x_4 + x_3 x_1}{x_1 x_2}.
$$

\subsection{The cluster category $\cat_{D_4}$}

The Auslander--Reiten quiver of $\cat_{D_4}$ is
$$
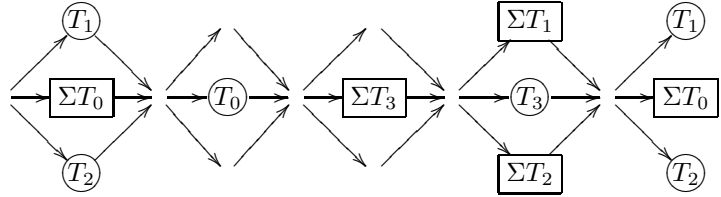
\xymatrix@-0.8pc{
  & *+[o][F]{T_1} \bdr & &  \bdr & & \bdr & &*+[F]{\susp T_1} \bdr & & *+[o][F]{T_1} \\
 \dr \hdr \bdr & *+[F]{\susp T_0} \dr & \dr \bdr \hdr & *+[o][F]{T_0} \dr &
\dr \bdr \hdr & *+[F]{\susp T_3}\dr & \dr \bdr \hdr & *+[o][F]{T_3} \dr &
 \dr \bdr \hdr &  *+[F]{\susp T_0} \\
 & *+[o][F]{T_2} \hdr & & \hdr & & \hdr & & *+[F]{\susp T_2} \hdr & & *+[o][F]{T_2} 
}
$$
The object $T:=T_1\oplus T_2 \oplus T_3 \oplus T_4$ is cluster-tilting.

The quiver of $B=\,$End$_{\cat_{D_4}}(T)$ is 
$$
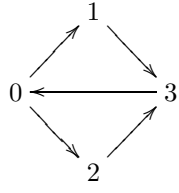
\xymatrix@-0.5pc{
  & 1 \bdr & \\
  0 \hdr \bdr & & 3 \ar@{->}[ll] \\
 & 2 \hdr &
}
$$
with the following relations: Any composition 
with the middle arrow vanishes, and the square
is commutative.

For $i=1,\ldots,n$, let $P_i$ be the image of $T_i$ in $\modb$, 
let $I_i$ be the image of $\susp^2 T_i$ and let $S_i$ be the simple 
top of $P_i$. Let $M$ and $N$ be the finite-dimensional $B$-modules given by: 
$$
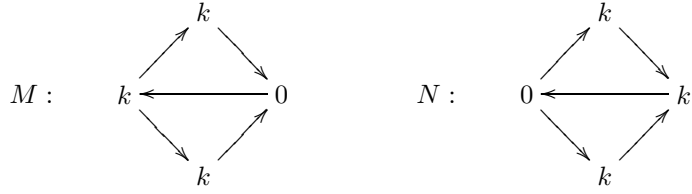
\xymatrix@-0.5pc{
& & k \bdr &  & & & & k \bdr & \\ 
M : & k \bdr \hdr & & 0 \ar@{->}[ll] & & N : & 0 \hdr \bdr & & k \ar@{->}[ll] \\
& & k \hdr & & & & & k \hdr & 
}
$$

As in the previous example, one can easily 
compute the AR-quiver of $B$.

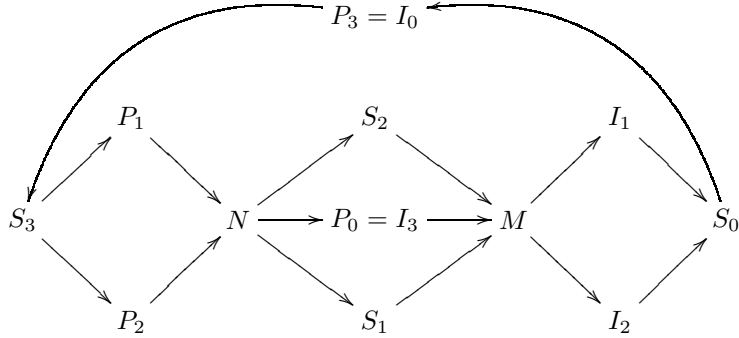
$$\xymatrix{
& & & P_3 = I_0 \ar@/_3pc/@{->}[ddlll] & & & \\
 & P_1 \bdr & & S_2 \bdr & & I_1 \bdr & \\
S_3 \bdr \hdr & & N \dr \bdr \hdr & P_0 = I_3 \dr & M \bdr \hdr & & S_0 \ar@/_3pc/@{->}[uulll] \\
& P_2 \hdr & & S_1 \hdr & & I_2 \hdr & 
}
$$
The submodules of $M$ are, up to isomorphism, 
$0$, $S_1$, $S_2$, $S_1\oplus S_2$ and $M$. 
Let $M_{\cat}$ be an indecomposable lift of $M$ in $\cat_{D_4}$.
Either by using add$\,T$-approximations and 
add$\,\susp T$-approximations or by~\cite[section 5.2]{KN}, 
one can compute the triangles
$$
T_3 \gfl T_0 \gfl M_\cat \gfl \susp T_3 
\;\; \text{ and } \;\; 
T_1\oplus T_2 \gfl T_0 \gfl \susm M \gfl \susp T_1\oplus \susp T_2.
$$ 
We thus have
$$
\ind M_\cat = [P_0]-[P_3] \text{, } \coi M_\cat = [P_1]+[P_2]-[P_0] 
$$
and
$$
X_{M_\cat} = \frac{(x_0 + x_3)^2  + x_1 x_2 x_3}{x_0 x_1 x_2}.
$$

\nocite{*}
\bibliographystyle{plain}
\bibliography{ClusterCharacters}

\end{document}